\documentclass[10pt,letterpaper]{amsart}

\usepackage[letterpaper, margin=1in]{geometry}
\usepackage{amssymb,amsthm,amsmath}
\usepackage{mathtools}
\usepackage[utf8]{inputenc}
\usepackage{enumerate}
\usepackage[T1]{fontenc}
\usepackage{tikz}
\usepackage{enumitem}
\usepackage{changepage}
\usepackage{tikz-cd}
\usepackage{textcomp}
\DeclareUnicodeCharacter{02BB}{\textquoteleft}
\usepackage{graphicx}
\usepackage[final]{microtype}
\usepackage[section]{placeins}
\usepackage{etoolbox}
\usepackage[numbers]{natbib}
\usepackage[pdftex,unicode]{hyperref}
\usepackage{lmodern}

\usetikzlibrary{arrows.meta, positioning}

\theoremstyle{plain}
\newtheorem{thm}[equation]{Theorem}
\newtheorem{lem}[equation]{Lemma}

\newtheorem{cor}[equation]{Corollary}

\theoremstyle{definition}
\newtheorem{definition}[equation]{Definition}
\newtheorem{example}[equation]{Example}

\theoremstyle{remark}
\newtheorem{rem}[equation]{Remark}

\newcommand{\ii}{\mathrm{i}}
\newcommand{\Z}{\mathbb{Z}}
\newcommand{\N}{\mathbb{N}}

\newcommand{\A}{\mathbb{A}}
\newcommand{\az}{\A_{\z}}
\newcommand{\q}{\mathbb{Q}}
\newcommand{\z}{\mathbb{Z}}
\newcommand{\ev}{\operatorname{ev}}

\newcommand{\F}{\mathbb{F}}

\newcommand{\p}{\mathbb{P}}

\newcommand{\bt}{\mathbf{t}}

\newcommand{\zhat}{\widehat{\Z}}
\newcommand{\nala}{\mathcal{L}_{\rm n}}
\newcommand{\arla}{\mathcal{L}_{\rm a}}
\newcommand{\adla}{\mathcal{L}}

\newcommand{\qhat}{\widehat{\q}}

\newcommand{\dov}{\mathrm{d}}

\newcommand{\ra}{\rightarrow}
\newcommand{\deff}{:=}

\newcommand{\OO}{\mathcal{O}}

\newcommand{\R}{\mathbb{R}}

\renewcommand{\a}{\mathrm{a}}

\renewcommand{\emph}[1]{{\it #1}}

\DeclareMathOperator{\Disc}{Disc}

\DeclareMathOperator{\id}{id}

\DeclareMathOperator{\End}{End}

\DeclareMathOperator{\Hom}{Hom}

\DeclareMathOperator{\td}{\mathrm{td}}

\DeclareMathOperator{\rk}{rk}
\DeclareMathOperator{\dima}{dim_a}
\DeclareMathOperator{\dimn}{dim_n}

\title[Finite-Dimensional Protori Are Adelic Tori]{Finite-Dimensional Protori Are Adelic Tori}
\author{Wayne Lewis}
\address{University of Hawaiʻi, Honolulu Community College}
\email{waynel@math.hawaii.edu}
\date{\today}


\begin{document}

\begin{abstract}
We show that every finite-dimensional compact connected abelian group is topologically isomorphic to an adelic torus $\frac{(\A^n/C)}{\Lambda}$ for the adeles $\A$, a torsion-free profinite group $C$, and a discrete cocompact $\Lambda\cong\q^n$.
A Lie theory is developed yielding a proper short exact sequence $0\ra \q^n \ra \mathcal{L}(G)\xrightarrow{\exp}G\ra 0$ for each adelic torus $G$ where $\exp$ denotes the adelic exponential and $\mathcal{L}(G)$ the adelic Lie algebra of $G$.
The adelic structure theory of the category $\mathbf{AT}$ of adelic tori under continuous quasi-homomorphism enables us to establish that the nonarchimedean dimension of an arbitrary number field $K$ with $[K\colon \q]=n$ is $1$, while the archimedean dimension is $n$.
\end{abstract}
\maketitle

\section{Introduction}
From its inception, the structure theory of compact connected abelian groups (also known as pro-tori \cite[Definition 9.30]{HofmannMorris}, herein contracted to protori) has been biased toward archimedean invariants.
Where previous topological investigations into protori---such as the resolution theorem and the characterisation of totally disconnected components---have largely treated nonarchimedean elements as secondary, the purpose of this paper is to establish a foundational, symmetric adelic topology for these spaces.

This adelic perspective clarifies the topological boundaries of these spaces and naturally connects their continuous geometry to the discrete arithmetic of number fields and the Krull-Schmidt category of discrete torsion-free abelian groups under quasi-homomorphism.

The adelic description is more than a change of coordinates.
It places the real and totally disconnected directions of a protorus in a common locally compact setting, where Pontryagin duality, self-dual resolutions, and an exponential sequence become aspects of the same structure. 
It also reveals a dimension that ordinary real dimension does not detect: the nonarchimedean dimension, measured by infinite procyclic factors of a profinite resolution. 
The distinction becomes particularly striking for number fields. 
A number field $K$ of degree $n$ has archimedean dimension $n$, yet its nonarchimedean dimension is always $1$, independently of $n$.
Thus adelic topology provides a common language for the continuous geometry of compact groups and the arithmetic of finite extensions of $\q$.

\subsection*{AI Declaration}
Anthropic Claude Opus 5.5 and OpenAI ChatGPT 6.2 Sol and Codex 26.1111 were used to assist with mathematical reasoning, Lean proof formalization and verification, and manuscript editing.
All mathematical claims remain the responsibility of the author.

\subsection*{Formal Verification}
Every result proved in this manuscript---specifically Theorems \ref{fdccg-adelic-tori}, \ref{spectrum-resolutions}, \ref{self-dual-pses}, \ref{symmetric-adelic-lie-algebra}, \ref{topological-characterisation-nonarchimedean-exponential}, \ref{adelic-exponential-lie-category}, \ref{kumu}, \ref{zassenhaus-realisation}, and \ref{nonarchimedean-line-number-field}, together with Lemma \ref{nala-periodic-lca}---has been formalized and machine-checked in Lean 4 with Mathlib (the only unformalized statement being Corollary \ref{FiveCatEquiv} which has proof following directly from cited results in \cite{DikranjanLewisLothMader,HofmannMorris}).
The formalization contains 234 audited declarations using only standard axioms.
The source is available at \url{https://github.com/WayneManAloha/AdelicTori} (commit \texttt{5401a88}).

\section{Structure of Adelic Tori}
Let $\mathbf{LCA}$ denote the category of locally compact abelian groups under continuous homomorphisms, and let $H^\vee$ denote the Pontryagin dual of $H\in\mathbf{LCA}$.
Define the {\it torsion-free rank} of a discrete finite-rank torsion-free abelian group $A$ by $\rk A\deff \dim_{\q}(\q\otimes_{\z}A)$; all tensor products herein are over $\z$.
Define the {\it archimedean dimension} of a protorus $G$ by $\dima G\deff \rk(G^\vee)$ \cite[Definition 8.23]{HofmannMorris}; a solenoid is a protorus of archimedean dimension $1$.

{\it Unless otherwise specified, all groups are abelian, and the symbol $\cong$ denotes topological isomorphism when it appears between topological groups.}

\begin{definition}
A \emph{discrete torus} is a discrete finite-rank torsion-free abelian group.
Discrete tori $A, B$ are \emph{quasi-isomorphic} if there is $n\in\N$, $f\colon A\ra B$, and $g\colon B\ra A$ such that $fg=n\cdot \id_B$ and $gf=n\cdot \id_A$.
\end{definition}

Define the elements of $\q\otimes\Hom(A,B)$ for discrete tori $A$ and $B$ to be the quasi-homomorphisms $A\ra B$; it is common practice in abelian group theory to abbreviate $\q\otimes A$ as $\q A$ and we make the identification $\q\otimes\Hom(A,B)=\q\Hom(A,B)$.
Discrete tori under quasi-homomorphism form an additive Krull-Schmidt category \cite[Theorem 7.4, Example 7.6]{Arnold}, which we denote $\mathbf{DT}$.
The category of finite-dimensional protori under continuous quasi-homomorphism, being contravariantly equivalent to $\mathbf{DT}$ via Pontryagin duality, is also an additive Krull-Schmidt category.

Let $G$ and $H$ be protori of archimedean dimension $n$.

\begin{definition}
An \emph{isogeny} from $G$ to $H$ is a surjective continuous homomorphism $G\ra H$ with finite kernel.
We say that $G$ and $H$ are \emph{isogenous}, and write $G\simeq H$, if there exists an isogeny $G\ra H$.
This forms an equivalence relation (symmetry and transitivity follow from taking pullbacks and pushouts of finite kernels).
There exists a continuous quasi-isomorphism $G\ra H$ if and only if there exists an isogeny $G\ra H$.
\end{definition}

A topological group $K$ is defined to be an $n$-dimensional {\it real torus}, $0\leq n\in\Z$, if and only if $K\cong \left(\frac{\R}{\z}\right)^n$ with the usual Euclidean topology on $\R$, subspace topology on $\Z$, quotient topology on $\frac{\R}{\z}$, and product topology on $\left(\frac{\R}{\z}\right)^n$.
A delta-subgroup $C$ of a protorus $G$ is a profinite subgroup such that $\frac{G}{C}$ is an $n$-dimensional real torus.

Define a finite-dimensional protorus $H$ to be {\it irreducible} if $H^\vee$ is strongly indecomposable; equivalently, $\q\End(H)$ is local \cite[Corollary 7.8]{Arnold}.

\begin{definition}
By \cite[Proposition 1]{Lewis2019}, every nonzero finitely generated profinite abelian group $D$ is topologically isomorphic to a unique \emph{standard representation} $\Delta = \prod_{j=1}^m \prod_{p\in\p}\zhat(p^{\mathfrak{n}_p(j)})$, where each $\Delta_j \deff \prod_{p\in\p}\zhat(p^{\mathfrak{n}_p(j)})$ ($1\le j\le m$) is a nontrivial procyclic group and $\Delta_{j+1}$ is a continuous quotient of $\Delta_j$.
The \emph{nonarchimedean dimension} of $D$ is $\dimn D \deff |\{j\in\{1,\dots,m\}\colon \Delta_j \text{ is infinite}\}|$, with $\dimn\{0\}\deff 0$ \cite[p.6]{Lewis2019}.
Define the {\it nonarchimedean dimension} of a finite-dimensional protorus $G$ to be $\dimn G \deff \dimn C$ for any delta-subgroup $C$ of $G$, which is well-defined and independent of the choice of $C$ \cite[p.8, Lemma 6 \& Corollary 4]{Lewis2019}.
\end{definition}

We set the {\it adelic dimension} of a finite-dimensional protorus $G$ to be the ordered pair $(\dimn(G),\dima(G))$, noting that $0\le\dimn(G)\le\dima(G)$ and that $\dimn(G)=0$ if and only if $\{0\}$ is a delta-subgroup of $G$, i.e., $G$ is an $n$-dimensional real torus \cite[Corollary 4]{Lewis2019}.

A {periodic $\mathbf{LCA}$ group $L$} is a totally disconnected $\mathbf{LCA}$ group equal to the directed union of its compact open subgroups; let $L_p$ denote the $p$-Sylow component of $L$ \cite[Definition 2.2]{HHR}.
A periodic $\mathbf{LCA}$ group $L$ with compact open subgroup $C$ is topologically isomorphic to the restricted product $\prod'_{p\in\p}(L_p,C_p)$ \cite[Theorem 3.3]{HHR}.
For a topological abelian group $L$ and a rational prime $p$, set $d_p(L)\deff\dim_{\q_p}\Hom(\q_p,L)$, where the continuous homomorphisms carry the natural $\q_p$-vector-space structure by precomposition with scalar multiplication on $\q_p$.
For a periodic $\mathbf{LCA}$ group $L$ that is not profinite, define its {\it nonarchimedean dimension} by
$\dimn L\deff\sup_{p\in\p}d_p(L)$.
The supremum is understood in the extended dimension sense if $d_p(L)$ is unbounded; in the finite-dimensional setting it is an integer.
For a profinite abelian group, the nonarchimedean dimension remains the dimension defined above by its standard procyclic representation; for delta-subgroups and protori, the preceding definitions remain unchanged.
Thus the two clauses apply to disjoint classes of periodic $\mathbf{LCA}$ groups.
The quantity $d_p(L)$ is invariant under topological isomorphisms, and so is the resulting dimension on nonprofinite periodic $\mathbf{LCA}$ groups.
For a finite-dimensional protorus $G$ of positive archimedean dimension, the periodic group $\nala(G)$ is nonprofinite, and the dimension formula below applies to it.
For the zero-dimensional case $G=\{0\}$, both Lie algebras have dimension zero.

Let $\zhat$ denote the additive topological group of the ring of profinite integers, defined as the inverse limit of the system $\left\{\frac{\z}{n\z}\colon 1\le n\in\Z\right\}$ directed by divisibility; then $\zhat\cong \prod_{p\in\p}\Z_p$ where $\Z_p$ denotes the ring of $p$-adic integers \cite[Example 2.3.11]{RZ}.
Let $\qhat = \prod'_{p\in\p}(\q_p,\z_p)$ denote the additive topological group of the locally compact ring of finite adeles; then $\qhat\cong \q\otimes\zhat$ \cite[Theorem 3.51]{HHR}.
Let $\A=\qhat\times\R$ denote the additive topological group of the locally compact ring of adeles; then $\A/\q\cong \q^\vee$ is the universal solenoid: it is the only torsion-free solenoid and every solenoid is topologically isomorphic to a quotient of $\A/\q$, whence a quotient of $\A$.

\begin{definition}
A topological group $G$ is an \emph{adelic torus} if $G\cong (\A^n/C)/\Lambda$ where $\Lambda\cong\q^n$ is embedded diagonally in $\A^n/C$ as a discrete cocompact subgroup for some $0\leq n\in\Z$ and a torsion-free profinite subgroup $C$ of $\A^n$.
Let $\mathbf{AT}$ denote the category of adelic tori under continuous quasi-homomorphism.
\end{definition}

\begin{thm}[Finite-Dimensional Compact Connected Abelian Groups Are Adelic Tori]\label{fdccg-adelic-tori}
If $G$ is a protorus with $\dima G=n$, then $G\cong (\A^n/C)/\Lambda$ for a torsion-free profinite group $C$ and a discrete cocompact $\Lambda\cong\q^n$.
\end{thm}
\begin{proof}
Let $G$ be a finite-dimensional protorus with $\dima G=n$.
Since $G^\vee$ is torsion-free of rank $n$, choose a full-rank free subgroup $F\cong\z^n$ of $G^\vee$ and identify $\z^n = F \subseteq G^\vee \subseteq \q^n$.
Endow $\A^n = \qhat^n \times \R^n$ with its standard Pontryagin self-duality pairing
\[
\langle (u,x), (v,y) \rangle = \prod_{p\in\p} e^{2\pi\ii \{u_p \cdot v_p\}_p} \cdot e^{-2\pi\ii \langle x, y \rangle},
\]
where $\{\cdot\}_p \colon \q_p \ra \q_p/\z_p \hookrightarrow \q/\z$ is the $p$-adic fractional part and $\langle x,y\rangle$ is the Euclidean inner product on $\R^n$.
Under this pairing, the diagonal subgroup $\widetilde{\Lambda} = \{(q,q)\colon q\in\q^n\} \subseteq \qhat^n\times\R^n = \A^n$ is discrete, cocompact, and satisfies $\widetilde{\Lambda}^\perp = \widetilde{\Lambda}$.

Viewing $G^\vee \subseteq \q^n \subseteq \qhat^n$ diagonally, let $C = \{u\in\qhat^n \colon \prod_{p\in\p} e^{2\pi\ii \{u_p\cdot b\}_p} = 1 \text{ for all } b\in G^\vee\}$ be the annihilator of $G^\vee$ in $\qhat^n$.
Because $\z^n \subseteq G^\vee$, $C$ is a closed subgroup of $(\z^n)^\perp \cap \qhat^n = \zhat^n$, so $C$ is a torsion-free profinite subgroup of $\qhat^n \subseteq \A^n$.
Identifying $C$ with $C\times\{0\}\subseteq\A^n$, we have $\widetilde{\Lambda} \cap C = \{0\}$.
Since $C$ is compact, the quotient map $\pi\colon \A^n \ra \A^n/C$ is closed with compact fibers, so $\Lambda \deff \pi(\widetilde{\Lambda}) \cong \q^n$ is discrete and cocompact in $\A^n/C$.

By Pontryagin duality, the dual of $(\A^n/C)/\Lambda$ is the annihilator of $C + \widetilde{\Lambda}$ in $(\A^n)^\vee \cong \A^n$, namely $(C^\perp \times \R^n) \cap \widetilde{\Lambda}^\perp = \{(q,q) \colon q\in\q^n \cap C^\perp\}$, where $C^\perp$ is the closure $\overline{G^\vee}$ of $G^\vee$ in $\qhat^n$.
For $q\in\q^n$, $q\in\overline{G^\vee}$ if and only if $q\in G^\vee + N\zhat^n$ for every $N\in\N$; since $\q^n \cap N\zhat^n = N\z^n \subseteq N G^\vee \subseteq G^\vee$, this holds if and only if $q\in G^\vee$.
Hence $\q^n \cap C^\perp = G^\vee$, so $\left((\A^n/C)/\Lambda\right)^\vee \cong G^\vee$, and Pontryagin reflexivity gives $G \cong (\A^n/C)/\Lambda$.
\end{proof}

Every abelian group $A$ decomposes as $A=\dov(A)\oplus M$ for the unique maximal divisible subgroup $\dov(A)$ and some complementary summand $M$ \cite[Theorem 2.5]{Fuchs}.
Fix $B\in \mathbf{DT}$.
Then $B$ is {\it locally free} if the localisation $B_{(p)}\cong\z_{(p)}\otimes_\z B$ (\cite[p.4]{Arnold}) is isomorphic to $\Z^{\rk B}_{(p)}$ for $p\in\p$; equivalently, $\dim_{\,\F_p}\left(\frac{B}{pB}\right)=\rk B$ for $p\in\p$ \cite[p.16]{Arnold}.
And $B$ is {\it quotient divisible} if there is $\Z^{\rk B}\cong F\subseteq B$ with ${\mathrm d}\left(\frac{B}{F}\right)=\frac{B}{F}$ \cite[p.18]{Arnold}.

The set of {\it supernatural numbers $\mathcal S$} is defined to be the set of formal products $\prod_{p\in\p}p^{k_p}$, $0\leq k_p\leq \infty$.
Define $\sim$ on $\mathcal S$ by $\bt\sim \bt'$ if and only if $k_p=\infty\Leftrightarrow k'_p=\infty$ and $k_p=k'_p$ for all but finitely many $p\in\p$; then $\sim$ is an equivalence relation on $\mathcal S$; the elements of $\frac{\mathcal S}{\sim}$, called (Baer) {\it types}, are partially ordered by $\bt\leq \bt'$ if and only if $k_p=\infty\Rightarrow k'_p=\infty$ and $k_p\leq k'_p$ for all but finitely many $p\in\p$; and one shows $\mathcal B\deff \frac{\mathcal S}{\sim}$ is a distributive lattice, called the Baer lattice.

Fix a rank-$1$ torsion free group $A$; fix a rational group $B$ (i.e., a group $B$ with $\Z\subseteq B\subseteq\q$) with $A\cong B$.
Define the type of a nonzero element $a\in A$ to be $\bt_A(a)\deff \frac{\prod_{p\in\p}p^{h_p}}{\sim}\in\mathcal B$ where $h_p\deff \sup\{\ell\in\N\cup\{0\}\colon a/p^\ell\in A\}$, $p\in\p$.
The (Baer) {\it type} of $A$ given by $\bt(A)\deff \bt_B(1)$ is well-defined.
A torsion-free group $E$ is $\bt(B)$-homogeneous if $\bt_E(a)=\bt(B)$ for all $0\neq a\in E$.

Set $\mathcal{D}(G)\deff\{C\colon C \text{ is a delta-subgroup of } G\}$ for $G\in\mathbf{AT}$.

\begin{thm}[Spectrum of Resolutions]\label{spectrum-resolutions}
Let $G$ be a finite-dimensional protorus, $\dim_a G=n$, and suppose $G\cong \frac{C'\times\R^n}{\z^n}$ for some delta-subgroup $C'\in \mathcal{D}(G)$.
If $\z^n\subseteq B\subseteq \q^n$, then $G\cong \frac{L_B\times\R^n}{B}$, where $L_B=C'\oplus_{\mathbb Z^n}B$ is the \textup{(}periodic $\mathbf{LCA}$\textup{)} pushout of $C'$ along $\mathbb Z^n\hookrightarrow B$ and $\Z^n\to C'$ is the projection of the diagonal embedding $\z^n\hookrightarrow C'\times\R^n$ onto its first coordinate.
\end{thm}
\begin{proof}
The existence of such a resolution associated with $B$ follows from the Structure Theorem for Finite-Dimensional Protori \cite[Theorem 1(3), $Y=B$]{Lewis2019}.
The following direct argument additionally identifies its periodic factor $L_B$ canonically, relative to the fixed resolution, as the stated pushout.

Let $\alpha\colon\z^n\ra C'$ be the projection onto the first coordinate of the diagonal embedding defining the given resolution, so that $G\cong \frac{C'\times\R^n}{\{(\alpha(z),z)\colon z\in\z^n\}}$.
Give $B$ the discrete topology and realize the pushout as $L_B=\frac{C'\oplus B}{\{(\alpha(z),-z)\colon z\in\z^n\}}$, with its quotient topology.
The subgroup $\{(\alpha(z),-z)\colon z\in\z^n\}$ is closed in $C'\oplus B$, since convergence in the discrete second coordinate forces $z$ to be eventually constant; hence, $L_B$ is an $\mathbf{LCA}$ group.
Write $\iota_B(b)=[0,b]\in L_B$, and embed $B$ diagonally in $L_B\times\R^n$ by $b\mapsto(\iota_B(b),b)$.

Define $\Phi_B\colon L_B\times\R^n\ra G$ by $\Phi_B([c,b],x)=[(c,x-b)]$.
The defining relation $[c+\alpha(z),b-z]=[c,b]$ shows that $\Phi_B$ is well-defined.
$\Phi_B$ is continuous and surjective and $\Phi_B(\iota_B(b),b)=0$ for every $b\in B$.
If $\Phi_B([c,b],x)=0$, then $(c,x-b)=(\alpha(z),z)$ for some $z\in\z^n$; hence, $x=b+z\in B$ and $[c,b]=[\alpha(z),b]=[0,b+z]=\iota_B(x)$; therefore, $\ker\Phi_B=\{(\iota_B(b),b)\colon b\in B\}$.
Thus, $\Phi_B$ induces a continuous isomorphism $\overline{\Phi}_B\colon\frac{L_B\times\R^n}{B}\ra G$.
To see directly that this is a topological isomorphism, define $C'\times\R^n\ra\frac{L_B\times\R^n}{B}$ by $(c,x)\mapsto [([c,0],x)]$.
For $z\in\z^n$ one has $[([\alpha(z),0],z)]=[([0,z],z)]=0$, so this map induces a continuous homomorphism $\Psi_B\colon G\ra\frac{L_B\times\R^n}{B}$.
The maps $\Psi_B$ and $\overline{\Phi}_B$ are inverse to one another.
It follows that $\overline{\Phi}_B$ is a topological isomorphism.

Finally, the map $L_B\ra \frac{B}{\z^n}$ by $[c,b]\mapsto b+\z^n$ gives an exact sequence $0\ra C'\ra L_B\ra \frac{B}{\z^n}\ra 0$.
Since $C'$ is profinite and $\frac{B}{\z^n}$ is a discrete torsion group, $L_B$ is a periodic $\mathbf{LCA}$ group; indeed, for each $\ell\in L_B$ some positive multiple $d\ell$ lies in $C'$, so $\overline{\langle\ell\rangle}$ is contained in finitely many translates of the profinite group $\overline{\langle d\ell\rangle}$.
\end{proof}

\section{Adelic Lie Theory}
Fix $A,B,C,D,E\in \mathbf{LCA}$.
A continuous homomorphism $f\colon A\ra B$ is {\it proper} if it is open onto its image.
A {\it proper exact sequence} in $\mathbf{LCA}$ is one for which all maps are proper.
A short exact sequence $0\ra C\ra D\ra E\ra0$ is a proper short exact sequence in $\mathbf{LCA}$ if and only if $0\ra E^\vee\ra D^\vee\ra C^\vee\ra0$ is a proper short exact sequence in $\mathbf{LCA}$ \cite[Proposition 2.2]{Moskowitz}.
Define a proper short exact sequence $0\ra C\ra D\ra E\ra0$ to be {\it Pontryagin self-dual} if $C\cong E^\vee$, $D\cong D^\vee$, and $E\cong C^\vee$ as topological groups.

\begin{thm}[Pontryagin Self-Dual Proper Short Exact Sequence]\label{self-dual-pses}
Let $G\in\mathbf{AT}$ with $\dima(G)=n$.
There is a periodic $\mathbf{LCA}$ group $L$ yielding a Pontryagin self-dual proper short exact sequence in $\mathbf{LCA}$\textup{:} $0\ra G^\vee\ra L\times\R^n\ra G\ra0$.
\end{thm}
\begin{proof}
Fix $\Delta\in\mathcal{D}(G)$ and a resolution $G \cong \frac{\Delta\times\R^n}{\{(\alpha(z),z)\colon z\in\z^n\}}$ where $\alpha\colon\z^n\ra\Delta$ has dense image.
Dualizing $0\ra\Delta\ra G\ra(\R/\z)^n\ra 0$ identifies $\z^n$ with a full-rank free subgroup $F \subseteq G^\vee$ satisfying $G^\vee/F \cong \Delta^\vee$.
Embedding $\z^n = F \subseteq G^\vee \subseteq \q^n$ and applying Theorem \ref{spectrum-resolutions} with $B = G^\vee$ gives the periodic $\mathbf{LCA}$ pushout $L = \Delta \oplus_F G^\vee$ and a short exact sequence $0\ra G^\vee\ra L\times\R^n\ra G\ra 0$.
Since $G^\vee$ is discrete and $L/\Delta \cong G^\vee/F$ is countable, $L\times\R^n$ is $\sigma$-compact; hence the continuous surjection $L\times\R^n\ra G$ onto the compact Hausdorff group $G$ is open \cite[Theorem 5.29]{HofmannMorris}, making the sequence proper in $\mathbf{LCA}$.

To see that $L \cong L^\vee$, note that $L \cong \prod'_{p\in\p}(L_p,\Delta_p)$ by \cite[Theorem 3.3]{HHR}.
For each $p\in\p$, $\alpha$ induces a continuous surjective $\z_p$-module homomorphism $\alpha_p\colon \z_p^n \twoheadrightarrow \Delta_p \cong \z_p^{s_p} \oplus \bigoplus_{i=1}^{r_p} (\z/p^{e_{p,i}}\z)$ ($e_{p,i}\ge 1$, $s_p + r_p \le n$), which is the Pontryagin dual of the $p$-primary embedding $(G^\vee/F)_p \hookrightarrow (\q_p/\z_p)^n$.
Because $\alpha_p$ is dual to $(G^\vee/F)_p \hookrightarrow (\q_p/\z_p)^n$, a single $\z_p$-basis of $\z_p^n$ adapts both $\alpha_p$ and $\z_p\otimes_\z G^\vee$, identifying $\z_p\otimes_\z G^\vee$ with $\q_p^{s_p} \oplus \bigoplus_{i=1}^{r_p} p^{-e_{p,i}}\z_p \oplus \z_p^{n-s_p-r_p}$ containing $\z_p^n$.
Pushing out $\Delta_p$ along $\z_p^n \hookrightarrow \z_p\otimes_\z G^\vee$ therefore yields $L_p \cong \q_p^{s_p} \oplus \bigoplus_{i=1}^{r_p} (\z/p^{2e_{p,i}}\z)$.
The standard Pontryagin self-duality isomorphism $L_p \xrightarrow{\sim} L_p^\vee$ carries the compact open subgroup $\Delta_p = \z_p^{s_p} \oplus \bigoplus_{i=1}^{r_p} p^{e_{p,i}}(\z/p^{2e_{p,i}}\z)$ onto its annihilator $\Delta_p^\perp \subseteq L_p^\vee$.
Hence these coordinatewise topological isomorphisms induce a topological isomorphism of restricted products $L^\vee \cong \prod'_{p\in\p}(L_p^\vee,\Delta_p^\perp) \cong \prod'_{p\in\p}(L_p,\Delta_p) \cong L$.
Since $(G^\vee)^\vee\cong G$ and $(\R^n)^\vee\cong\R^n$, the proper short exact sequence is Pontryagin self-dual.
\end{proof}

Let $H$ be a Hausdorff topological group.
Define the {\it totally disconnected component of $0$ in $H$}, denoted $\td(H)$, to be the set of {\it quasi-torsion elements} $x\in H$: those $x\in H$ such that $\langle x\rangle$ is either finite or its subspace topology is non-discrete and linear.
If $H$ is also complete, as when $H\in\mathbf{LCA}$, then $\td(H)$ coincides with the union of all profinite subgroups of $H$ \cite[Proposition 9(3)]{DikranjanLewisLothMader}.

Fix $G\in\mathbf{AT}$ with $\dima(G)=n$.
Define the {\it nonarchimedean Lie algebra of $G$}, denoted $\nala(G)$, to be $\Hom(\zhat,G)$ under the compact-open topology.
We prove in Lemma \ref{nala-periodic-lca} below that $\nala(G)$ is a periodic $\mathbf{LCA}$ group.
Lemma \ref{nala-periodic-lca} shows that evaluation at $1$ induces an algebraic isomorphism $\nala(G)\cong\td(G)$; thus, $\nala(G)$ is a divisible periodic $\mathbf{LCA}$ group uniquely determined up to topological isomorphism by its restricted direct product decomposition \cite[Theorem 3.3]{HHR}.
For every $p\in\p$, the continuous-homomorphism and Pontryagin-duality argument in Theorem \ref{symmetric-adelic-lie-algebra} gives $d_p(\nala(G))=n$, hence $\dimn\nala(G)=n$.
Define the {\it archimedean Lie algebra of $G$}, denoted $\arla(G)$, to be $\Hom(\R,G)$ under the compact-open topology; then $\arla(G)\cong\R^{\dima(G)}$ as locally compact real topological vector spaces \cite[Theorem 8.22(6)]{HofmannMorris}.
Thus, $(\dimn\nala(G),\dima\arla(G))=(n,n)$, whereas the adelic dimension of $G$ itself is $(\dimn G,n)$ and often has $\dimn G<n$; for an $n$-dimensional real torus, $\dimn G=0$ while $\dimn\nala(G)=n$.
Set $\az\deff\zhat\times\R$.
Define the {\it adelic Lie algebra of $G$}, denoted $\adla(G)$, to be $\Hom(\az,G)$ under the compact-open topology.

A topological group $G$ {\it has no small subgroups} if there is a neighborhood $U$ of $0$ such that for every subgroup $H\subseteq U$ one has $H=\{0\}$; for example, $\R^n$ and $\left(\frac{\R}{\z}\right)^n$ have no small subgroups for $n\in\N$.

\begin{lem}\label{nala-periodic-lca}
Let $G\in\mathbf{AT}$ and endow $\Hom(\zhat,G)$ with the compact-open topology.
Then $\Hom(\zhat,G)$ is a periodic $\mathrm{LCA}$ group, and evaluation at $1$ induces an algebraic isomorphism $\ev_1\colon\Hom(\zhat,G)\ra\td(G)$.
\end{lem}
\begin{proof}
Fix $\Delta\in\mathcal D(G)$ and set $K_\Delta=\{f\in\Hom(\zhat,G)\colon f(\zhat)\subseteq\Delta\}$.
We prove that $K_\Delta$ is a compact open subgroup of $\Hom(\zhat,G)$.

First, $K_\Delta$ is open:
Let $\pi_\Delta\colon G\ra\frac{G}{\Delta}$ be the quotient map.
Since $\Delta\in\mathcal D(G)$, one has $\frac{G}{\Delta}\cong\left(\frac{\R}{\z}\right)^n$, $n=\dima(G)$, and therefore $\frac{G}{\Delta}$ has no small subgroups.
Choose an identity neighborhood $U\subseteq\frac{G}{\Delta}$ containing no nontrivial subgroup.
Then $\left\{\varphi\in\Hom\left(\zhat,\frac{G}{\Delta}\right)\colon\varphi(\zhat)\subseteq U\right\}=\{0\}$.
This is an open neighborhood of $0$ in the compact-open topology, so $\Hom\left(\zhat,\frac{G}{\Delta}\right)$ is discrete.
Define $\rho_\Delta\colon\Hom(\zhat,G)\ra\Hom\left(\zhat,\frac{G}{\Delta}\right)$ by $\rho_\Delta(f)=\pi_\Delta\circ f$.
The map $\rho_\Delta$ is continuous and $K_\Delta=\ker(\rho_\Delta)$.
Since the codomain is discrete, $K_\Delta$ is open.

We next prove that $K_\Delta$ is compact:
Since $\Delta$ is profinite, its canonical continuous $\zhat$-action shows that evaluation at $1$ induces a topological isomorphism $\ev_1\colon K_\Delta=\Hom(\zhat,\Delta)\ra\Delta$.
Indeed, $\ev_1$ is injective because $\zhat$ is topologically generated by $1$, while its inverse sends $x\in\Delta$ to the continuous homomorphism $z\mapsto z\cdot x$.
Thus, $K_\Delta\cong\Delta$ is a compact open subgroup of $\Hom(\zhat,G)$.

It remains to prove periodicity and identify evaluation at $1$.
By \cite[Proposition 9(3)]{DikranjanLewisLothMader}, $\td(G)$ is the union of all profinite subgroups of $G$; for any profinite subgroup $P$ of $G$, $P+\Delta$ is profinite and $G/(P+\Delta)$ is a quotient of $G/\Delta\cong(\R/\z)^n$ by the finite subgroup $(P+\Delta)/\Delta$, hence an $n$-dimensional real torus, so $P\subseteq P+\Delta\in\mathcal D(G)$ and $\td(G)=\bigcup_{\Delta'\in\mathcal D(G)}\Delta'$.
Let $f\in\Hom(\zhat,G)$.
The homomorphic image $f(\zhat)$ is profinite, so $f(1)\in\td(G)$, and we may choose $\Delta_f\in\mathcal D(G)$ such that $f(1)\in\Delta_f$.
Since $f(m)=m f(1)$ for $m\in\Z$ and $\Z$ is dense in $\zhat$, continuity gives $f(z)=z\cdot f(1)\in\Delta_f$ for every $z\in\zhat$.
Hence, $f\in K_{\Delta_f}$.
Thus, every element of $\Hom(\zhat,G)$ belongs to a compact subgroup, so $\Hom(\zhat,G)$ is a periodic $\mathbf{LCA}$ group.

Evaluation at $1$ is injective because $\zhat$ is topologically generated by $1$.
Conversely, if $x\in\td(G)$, then $x$ lies in some $\Delta'\in\mathcal D(G)$, and the canonical $\zhat$-action on $\Delta'$ defines $f_x\colon\zhat\ra G$ by $f_x(z)=z\cdot x$.
Then $f_x(1)=x$.
Thus, evaluation at $1$ is an algebraic isomorphism onto $\td(G)$.
\end{proof}

Our use of nonarchimedean scalar notation depends fundamentally on \cite[Proposition 4.21]{HHR}: if $L$ is a periodic $\mathbf{LCA}$ group, then the closure of the cyclic subgroup generated by each $g\in L$ is procyclic and has the Sylow decomposition $\overline{\langle g\rangle}\cong\prod_{p\in\p}\overline{\langle g\rangle}_p$.
Consequently, the Sylow decomposition induces a canonical continuous bilinear action $\zhat\times L\ra L$ via $(z,g)\mapsto z\cdot g$, where, for $z=(z_p)_{p\in\p}\in\zhat=\prod_{p\in\p}\z_p$, the $p$-component of $z\cdot g$ is $(z\cdot g)_p=z_pg_p$.
Thus, the scalar notation used below is not an additional multiplication imposed on the additive group $L$, but its canonical external $\zhat$-action.
In particular, delta-subgroups of adelic tori and compact open subgroups of $\nala(G)$ are compact topological $\zhat$-modules.

\begin{rem}[Archimedean and Nonarchimedean Scalar Actions]
The canonical $\zhat$-action on a periodic $\mathbf{LCA}$ group is the nonarchimedean analogue of the familiar scalar action on the archimedean Lie algebra and its non-locally compact copy ${\mathrm a}(G)$ in $G$.
In the representation of a finite-dimensional protorus given by the Resolution Theorem \cite[Theorem 8.20]{HofmannMorris}, the archimedean component $\arla(G)\cong\R^n$ carries its natural $\R$-module structure \cite[Theorem 8.22(6)]{HofmannMorris}.
Accordingly, we use $\R$-scalars on archimedean components and $\zhat$-scalars on periodic nonarchimedean components.
On torsion-free divisible periodic $\mathbf{LCA}$ groups, the restricted direct product decomposition extends the canonical $\zhat$-action coordinatewise to a canonical $\qhat$-action.
No additional multiplicative structure is imposed on the ambient $G\in\mathbf{AT}$.
\end{rem}

\begin{thm}[Symmetric Adelic Lie Algebra]\label{symmetric-adelic-lie-algebra}
Let $G\in\mathbf{AT}$ with $\dima(G)=n$.
The adelic Lie algebra of $G$ decomposes algebrotopologically in $\mathbf{LCA}$ as
$\adla(G)\cong\nala(G)\times\arla(G)$.
Here $\nala(G)$ is a divisible periodic $\mathbf{LCA}$ group with $d_p(\nala(G))=n$ for every $p\in\p$, so $\dimn\nala(G)=n$; also, $\arla(G)\cong\R^n$ as a real topological vector space.
The restricted-product decomposition of the periodic $\mathbf{LCA}$ group $\nala(G)$ is given by \cite[Theorem 3.3]{HHR}.
\end{thm}
\begin{proof}
Using the canonical decomposition $\az=\zhat\times\R$, define $\Phi\colon\Hom(\az,G)\ra\Hom(\zhat,G)\times\Hom(\R,G)$ by $\Phi(f)=(f\vert_{\zhat},f\vert_{\R})$.
Its inverse is given by $\Psi(u,v)(z,x)=u(z)+v(x)$.
Thus, $\Phi$ is an algebraic isomorphism.

The restriction maps defining $\Phi$ are continuous in the compact-open topology.
To prove continuity of $\Psi$, let $K\subseteq\az$ be compact and let $U$ be an identity neighborhood in $G$.
Choose identity neighborhoods $U_{\mathrm n},U_{\mathrm a}\subseteq G$ such that $U_{\mathrm n}+U_{\mathrm a}\subseteq U$.
Let $K_{\mathrm n}\subseteq\zhat$ and $K_{\mathrm a}\subseteq\R$ be the projections of $K$.
Then $K_{\mathrm n}$ and $K_{\mathrm a}$ are compact and
$\Psi([K_{\mathrm n},U_{\mathrm n}]\times[K_{\mathrm a},U_{\mathrm a}])\subseteq[K,U]$.
Hence, $\Psi$ is continuous; therefore, $\adla(G)\cong\nala(G)\times\arla(G)$.

By Lemma \ref{nala-periodic-lca}, evaluation at $1$ induces an algebraic isomorphism $\nala(G)\cong\td(G)$ and $\nala(G)$ is a divisible periodic $\mathrm{LCA}$ group and for any $\Delta\in\mathcal D(G)$ the compact open subgroup $K_\Delta$ of $\nala(G)$ is topologically isomorphic to $\Delta$.
For each $p\in\p$, evaluation at $1\in\zhat$ yields an isomorphism of $\q_p$-vector spaces
$\Hom(\q_p,\nala(G))\cong\Hom(\q_p,G)$.
Indeed, the forward map sends $u$ to $x\mapsto u(x)(1)$.
In the reverse direction, a continuous homomorphism $v\colon\q_p\ra G$ determines $x\mapsto[z\mapsto v(z_p x)]$, where $z_p$ denotes the $p$-component of $z\in\zhat$; continuity follows from the compact-open exponential law.
The two maps are inverse because $1$ topologically generates $\zhat$.

Set $A=G^\vee$, a discrete torsion-free group of rank $n$.
The self-duality of the additive locally compact group $\q_p$ and Pontryagin duality identify the foregoing vector space with $\Hom(\q_p,G)\cong\Hom(A,\q_p)\cong \Hom_{\q}(\q\otimes_{\z}A,\q_p)\cong\q_p^n$.
In the middle step, a homomorphism from $A$ to the divisible group $\q_p$ extends uniquely to its rational vector-space span; all maps out of the discrete group $A$ are continuous.
It follows for every prime $p$ that $d_p(\nala(G))=\dim_{\q_p}\Hom(\q_p,\nala(G))=n$; hence, $\dimn\nala(G)=\sup_p d_p(\nala(G))=n$.
\end{proof}

\begin{example}
For the real torus $G=(\R/\z)^n$, one has $\dimn G=0$ but
$\nala(G)\cong(\q/\z)^n$ with the discrete compact-open topology.
For every prime $p$, the theorem gives
$\dim_{\q_p}\Hom(\q_p,(\q/\z)^n)=n$.
Thus $\dimn\nala(G)=n$ even though $\dimn G=0$.
\end{example}

Fix $G\in\mathbf{AT}$ with $\dima(G)=n$.
Define the {\it nonarchimedean exponential} $\exp_{\mathrm n}\colon\nala(G)\ra G$ by $\exp_{\mathrm n}(f)=f(1)$.
Define the {\it archimedean exponential} $\exp_{\mathrm a}\colon\arla(G)\ra G$ by $\exp_{\mathrm a}(g)=g(1)$.
Define the {\it adelic exponential} $\exp\colon\adla(G)\ra G$ by $\exp(f,g)=f(1)+g(1)$.

The {\it arc component of $0$ in $G\in\mathbf{AT}$}, denoted ${\mathrm a}(G)$, is the canonical fully invariant dense divisible arc-connected component of $0$ in $G$ \cite[Corollary 8.47]{HofmannMorris}.

A rank-$n$ group $B\in\mathbf{DT}$ is {\it strongly indecomposable} if $B$ is not quasi-isomorphic to a direct sum of nontrivial discrete tori of ranks less than $n$; equivalently, $\q\End(B)$ is local \cite[Corollary 7.8]{Arnold}.
We assume that quasi-homomorphism groups and quasi-endomorphism rings of adelic tori consist only of continuous quasi-homomorphisms and are endowed with the compact-open topology.

For an abelian group $A$ and $p\in\p$, one has $\frac{A}{pA}\cong\frac{\z}{p\z}\otimes_\z A$, and we define
$p\text{-rank}(A)\deff\dim_{\z/p\z}\frac{A}{pA}$.
The delta-subgroups of an irreducible adelic torus $G$ of adelic dimension $(m,n)$ are commensurable and form a semilattice with directed union $\td(G)$ \cite[Theorem 6(1) \& Proposition 9(3)]{DikranjanLewisLothMader}.

\begin{thm}[Topological Characterisation via Nonarchimedean Exponential]\label{topological-characterisation-nonarchimedean-exponential}
The nonarchimedean exponential is a continuous algebraic isomorphism of $\nala(G)$ onto $\td(G)$.
If $G\neq\{0\}$, then $\td(G)$ is not locally compact in its subspace topology and $\exp_{\mathrm n}\colon\nala(G)\ra\td(G)$ is not open.
For each $\Delta\in\mathcal D(G)$, the nonarchimedean exponential maps the compact open subgroup
$K_\Delta=\{f\in\nala(G)\colon f(\zhat)\subseteq\Delta\}$ topologically onto the compact non-open subgroup $\Delta$ of $\td(G)$.
\end{thm}
\begin{proof}
For $f\in\nala(G)$, the homomorphic image $f(\zhat)$ is profinite, and hence $f(1)\in\td(G)$.
Thus, $\exp_{\mathrm n}$ is well-defined as a map into $\td(G)$.
Evaluation at $1$ is continuous in the compact-open topology.
It is injective because $\zhat$ is topologically generated by $1$.

Let $x\in\td(G)$.
Then $x$ lies in a delta-subgroup $\Delta\in\mathcal D(G)$.
The canonical continuous $\zhat$-action on $\Delta$ defines a continuous homomorphism $f_x\colon\zhat\ra G$ by $f_x(z)=z\cdot x$.
Since $\exp_{\mathrm n}(f_x)=f_x(1)=x$, the map $\exp_{\mathrm n}$ is surjective onto $\td(G)$.
Thus, $\exp_{\mathrm n}$ is a continuous algebraic isomorphism of $\nala(G)$ onto $\td(G)$.

For $\Delta\in\mathcal D(G)$, Lemma \ref{nala-periodic-lca} shows that $K_\Delta$ is compact and open in $\nala(G)$ and that evaluation at $1$ restricts to a topological isomorphism $K_\Delta\ra\Delta$.

Suppose that $\Delta$ is open in the subspace topology on $\td(G)$.
Then there is an identity neighborhood $U$ in $G$ such that $U\cap\td(G)\subseteq\Delta$.
Since $\td(G)$ is dense in $G$, the set $U\cap\td(G)$ is dense in $U$.
Since $\Delta$ is compact and therefore closed in $G$, it follows that $U\subseteq\Delta$.
Thus, $\Delta$ is open in $G$.
As $G$ is connected, one obtains $\Delta=G$, which is impossible when $G\neq\{0\}$ because $\Delta$ is totally disconnected.
Therefore, $\Delta$ is not open in $\td(G)$, and consequently $\exp_{\mathrm n}$ is not open.

Finally, suppose that $\td(G)$ is locally compact in its subspace topology.
Choose an identity neighborhood $W$ in $\td(G)$ having compact closure $C$ in $\td(G)$.
There is an identity neighborhood $U$ in $G$ such that $U\cap\td(G)\subseteq W\subseteq C$.
The subgroup $\td(G)$ is dense in $G$, while $C$ is compact and therefore closed in $G$.
Hence, $U\subseteq C\subseteq\td(G)$.
Thus, $\td(G)$ is an open subgroup of the connected group $G$, so $\td(G)=G$.
Since $\td(G)$ is totally disconnected, this forces $G=\{0\}$.
Therefore, $\td(G)$ is not locally compact when $G\neq\{0\}$.
\end{proof}

\begin{thm}[Adelic Exponential Yields Lie Category]
\label{adelic-exponential-lie-category}
Under the decomposition
$\adla(G)\cong\nala(G)\times\arla(G)$,
the adelic exponential is the resolution map of
\cite[Theorem 15]{DikranjanLewisLothMader},
whence continuous, surjective, and open with
discrete kernel isomorphic to $\q^n$.
\end{thm}
\begin{proof}
This follows from \cite[Theorem 15]{DikranjanLewisLothMader}.
\end{proof}

Let $G\in\mathbf{AT}$.
The {\it kumu} of $G$ is defined to be $\kappa(G)\deff\td(G)\cap{\mathrm a}(G)$.
$G$ is {\it torus-free} if it contains no nontrivial subgroups topologically isomorphic to a torus.

\begin{thm}[The Kumu]\label{kumu}
The kumu $\kappa(G)$ of a torus-free $G\in\mathbf{AT}$ is the canonical fully invariant dense subgroup of $G$ algebraically isomorphic to $\q^{\dima G}$.
\end{thm}
\begin{proof}
This follows from \cite[Theorem 15]{DikranjanLewisLothMader}.
\end{proof}

Fix $G\in\mathbf{AT}$ with $\dima(G)=n$.
Let $B\subseteq\kappa(G)$ be dense in $G$, and suppose that $A\in\mathbf{DT}$ is isomorphic to $B$.
Fix an isomorphism $\alpha\colon A\ra B$.
The {\it dual topology on $A$ relative to $(G,\alpha)$} is the topology transported to $A$ by $\alpha$, where $B$ carries the subspace topology inherited from $G$.
When no reference to ``relative to $(G,\alpha)$'' is made, we adopt the convention that for a realization $A\subseteq\q^n\subseteq\R^n$ with no free summands, identifying $\arla(A^\vee)=\Hom(\R,A^\vee)$ with $\R^n$ via the standard Euclidean inner product $\langle x,a\rangle \mapsto e^{2\pi\ii\langle x,a\rangle}$ ($x\in\R^n$, $a\in A$), the dual topology on $A$ is its dual topology relative to $(A^\vee,\exp_{\mathrm a}\vert_A)$.

\begin{rem}[Dense Subgroups and Dual Topologies]\label{dense-subgroups-dual-topologies}
Let $A\in\mathbf{DT}$ have rank $n\in\N$ and no free summands, and let $H$ be the closure of $A$ in $\R^n$.
By the structure theorem for closed subgroups of $\R^n$, $H=V\oplus\Lambda$ for a real vector subspace $V$ and a discrete free subgroup $\Lambda$; projecting $A$ onto $\Lambda$ shows that $\Lambda=\{0\}$ \textup{(}as $A$ has no free summands\textup{)}, and since $\q A=\q^n$, we have $H=\R^n$.
Because $A$ has no free summands, $\Hom(A,\z)=0$, so $\exp_{\mathrm a}\colon\R^n\ra A^\vee$ is injective \textup{(}cf.~\cite[Corollary 8.47]{HofmannMorris}\textup{)} and restricts to an isomorphism $B\deff\exp_{\mathrm a}(A)\cong A$.
Since $\exp_{\mathrm a}$ is continuous with dense image ${\mathrm a}(A^\vee)$ \cite[Corollary 8.47]{HofmannMorris} and satisfies $\exp_{\mathrm a}(\q^n)=\kappa(A^\vee)$ \cite[Theorem 15]{DikranjanLewisLothMader}, $B$ is dense in the compact Hausdorff \textup{(}hence complete\textup{)} group $A^\vee$, so the completion of $A$ under its dual topology relative to $(A^\vee,\exp_{\mathrm a}\vert_A)$ is topologically isomorphic to $A^\vee$.
\end{rem}

\begin{cor}[Five Equivalents of Topological Isomorphism of Adelic Tori]\label{FiveCatEquiv}
TFAE for torus-free $G, H\in \mathbf{AT}$:
\begin{enumerate}[label=(\roman*), font=\upshape]
\item $G\cong H$;
\item $\td(G)\cong\td(H)$;
\item $\kappa(G)\cong\kappa(H)$;
\item $\a(G)\cong\a(H)$;
\item $G$ is homeomorphic to $H$.
\end{enumerate}
\end{cor}
\begin{proof}
(i)$\Leftrightarrow$(ii): \cite[Theorem 12]{DikranjanLewisLothMader}.\newline
(ii)$\Rightarrow$(iii): $\td(G)\cong\td(H)\Rightarrow G\cong H$ so $\kappa(G)\cong \kappa(H)$ by Theorem \ref{kumu}.\newline
(iii)$\Rightarrow$(iv): $\kappa(G)$ and $\kappa(H)$ are dense in $G$ and $H$, respectively, by Theorem \ref{kumu}.
Let $f\colon \kappa(G)\ra\kappa(H)$ be a topological isomorphism.
Because $H$ is compact and therefore complete, and $\kappa(G)$ is dense in $G$, there is a unique continuous extension $\bar f\colon G \ra H$.
Symmetric unique extension of $f^{-1}$ yields continuous inverses, so $G\cong H$.
It follows that $\a(G)\cong \a(H)$.\newline
(iv)$\Rightarrow$(i): Let $f\colon \a(G)\ra\a(H)$ be a topological isomorphism.
The group $H$ is compact, hence complete, and $\a(G)$ is dense in $G$.
Hence, there is a unique extension morphism $\bar f\colon G \ra H$ of $f$.
Similarly, we have the unique continuous extension $\overline {f^{-1}}\colon H\ra G$ of $f^{-1}$.
The morphism $\overline {f^{-1}} \circ \bar f\colon G \ra G$ extends $\mathrm{id}_{\a(G)}\colon \a(G)\ra \a(G)$, which is also extended by $\mathrm{id}_G$; hence, by uniqueness, we have $\overline {f^{-1}}\circ \bar f = \mathrm{id}_G$.
Similarly, $\bar f \circ \overline {f^{-1}} = \mathrm{id}_H$.
Hence, $(\bar f)^{-1} = \overline {f^{-1} }$ is continuous.
It follows $G\cong H$, whence $G$ is homeomorphic to $H$.\newline
(v)$\Leftrightarrow$(i) \cite[Theorem 8.58]{HofmannMorris}.
\end{proof}

\section{Adelic Lines and Number Fields}
Let $A\in\mathbf{DT}$ have rank $n$ and no free summands.
Let $F\cong\z^n$ be a full-rank free subgroup of $A$ and set $T\deff\frac{A}{F}$.
Then there is a proper short exact sequence $0\ra F\ra A\ra T\ra0$, whose Pontryagin dual is the proper short exact sequence $0\ra T^\vee\ra A^\vee\ra\left(\frac{\R}{\z}\right)^n\ra 0$.
Since $\frac{A^\vee}{T^\vee}\cong F^\vee\cong\left(\frac{\R}{\z}\right)^n$, the profinite group $T^\vee$ is a delta-subgroup of $A^\vee$.
In contrast to a resolution supported on $T^\vee$, $\nala(A^\vee)\times\R^n$ is the divisible locally compact cover associated with the adelic exponential, whose diagonal kernel is isomorphic to $\q^n$, and algebraically $\nala(A^\vee)\cong\td(A^\vee)$.

Now, fix a number field $K$, that is, a finite extension field of $\q$.
A {\it group of integers of $K$} is a torsion-free subgroup $B_K$ with no free summands satisfying $\OO_K\subseteq B_K\subseteq K$, $\q B_K=K$, and $\End(B_K)=\OO_K$ \cite[Theorem, p.180]{Zassenhaus}; the terminology is endomorphism-theoretic: the word ``integers'' refers to the endomorphisms of $B_K$, not to the assertion that all elements of $B_K$ are algebraic integers.
Let $\R\cdot K$ denote Zassenhaus' {\it hypercomplex system with the same basis and multiplication constants as $\OO_K$ over the rational number field $\q$} \cite[Theorem, p.180]{Zassenhaus}.
After choosing a $\z$-basis of $\OO_K$, we identify $\OO_K^\vee$ topologically with $(\R\cdot K)/\OO_K$.
Set $T_K\deff B_K/\OO_K$.

The next theorem is an adelic realisation of Zassenhaus' construction \cite[Theorem, p.180]{Zassenhaus}, adapted so that each stage has prime index $q_k$ and cyclic quotient $\z/q_k\z$ for an arbitrary number field $K$.

\begin{thm}[Adelic Realisation of a Group of Integers]\label{zassenhaus-realisation}
Every number field $K$, with $n=[K\colon\q]$, contains a homogeneous, strongly indecomposable, locally free torsion-free subgroup $B_K = \sum_{k=1}^\infty \mathfrak{p}_k^{-1}$ \textup{(}homogeneous of type $\bt(\z)$ when $n\ge 2$\textup{)} such that $\OO_K\subseteq B_K\subseteq K$, $\q B_K=K$, $\End(B_K)=\OO_K$, and $(B_K/\OO_K)_p$ is cyclic of prime order or trivial for every $p\in\p$.

Moreover, setting $T_K\deff B_K/\OO_K$, $T_K^\vee$ is a delta-subgroup of $B_K^\vee\in\mathbf{AT}$, and the projection of the diagonal copy of $\OO_K$ into $T_K^\vee$ is dense.
After choosing a $\z$-basis of $\OO_K$, there are topological isomorphisms $B_K^\vee\cong \frac{T_K^\vee\times\R\cdot K}{\OO_K}\cong\frac{\nala(B_K^\vee)\times\R\cdot K}{K}$, and $\frac{B_K^\vee}{T_K^\vee}\cong\frac{\R\cdot K}{\OO_K}\cong\left(\frac{\R}{\z}\right)^{[K\colon\q]}$.
\end{thm}
\begin{proof}
If $n=1$, then $K=\q$ and $\OO_\q=\z$; choosing an infinite proper subset $\{q_1,q_2,\dots\}\subsetneq\p$ and setting $\mathfrak{p}_k = q_k\z$ gives $B_\q = \sum_{k=1}^\infty q_k^{-1}\z \subsetneq \q$, which is locally free, homogeneous, strongly indecomposable, has no free summands \textup{(}as $\rk B_\q=1$ and $B_\q\not\cong\z$\textup{)}, and satisfies $\End(B_\q) = \{x\in\q\colon x B_\q\subseteq B_\q\} = \z$ \textup{(}since no prime has infinite height in $B_\q$\textup{)} and $B_\q/\z \cong \bigoplus_{k=1}^\infty \z/q_k\z$.

Assume now that $n=[K\colon\q]\ge 2$.
Let $L/\q$ be a Galois closure of $K/\q$, and let $\sigma_1,\dots,\sigma_n \colon K \hookrightarrow L$ be the $n$ distinct field embeddings of $K$ into $L$.
Choose a primitive element $t\in\OO_K$ ($K=\q(t)$) with monic minimal polynomial $P(x)\in\z[x]$ of degree $n$.
The $L$-algebra map $\Phi\colon L\otimes_{\q} K \ra L^n$ given by $a\otimes b \mapsto (a\sigma_1(b),\dots,a\sigma_n(b))$ is an isomorphism because $\sigma_1(t),\dots,\sigma_n(t)$ are distinct.

Since infinitely many rational primes split completely in $L$ \cite[Chapter 8, Theorem 43]{Marcus}, the set of rational primes $q$ that split completely in $L$ and do not divide $[\OO_K \colon \z[t]]$ is infinite; enumerate this set as a sequence $\{q_1,q_2,\dots\}$ of pairwise distinct primes.
For each $k\ge 1$, fix a prime ideal $\mathfrak{P}_k$ of $\OO_L$ above $q_k$, so $\OO_L/\mathfrak{P}_k \cong \F_{q_k}$.
Because $q_k$ splits completely in $L$, $q_k \nmid \Disc(\OO_K)$, and since $q_k \nmid [\OO_K \colon \z[t]]$, we have $q_k \nmid \Disc(P)$, so $\overline{\sigma_1(t)},\dots,\overline{\sigma_n(t)} \in \OO_L/\mathfrak{P}_k \cong \F_{q_k}$ are pairwise distinct.
Let $i(k)\in\{1,\dots,n\}$ satisfy $k\equiv i(k)\pmod n$, and set $\mathfrak{p}_{k,j} \deff \OO_K \cap \sigma_j^{-1}(\mathfrak{P}_k)$ for $1\le j\le n$ and $\mathfrak{p}_k \deff \mathfrak{p}_{k,i(k)}$.
Then $\mathfrak{p}_{k,1},\dots,\mathfrak{p}_{k,n}$ are the $n$ distinct degree-$1$ prime ideals of $\OO_K$ above $q_k$, so $q_k\OO_K = \prod_{j=1}^n \mathfrak{p}_{k,j}$, $\mathfrak{p}_k^{-1}/\OO_K \cong \OO_K/\mathfrak{p}_k \cong \z/q_k\z$, and the Chinese Remainder Theorem gives an isomorphism
\[
\pi_k \colon \frac{\OO_K}{q_k\OO_K} \xrightarrow{\sim} \prod_{j=1}^n \frac{\OO_K}{\mathfrak{p}_{k,j}} \xrightarrow{\sim} \F_{q_k}^n, \qquad b + q_k\OO_K \mapsto (\overline{\sigma_1(b)},\dots,\overline{\sigma_n(b)}) \pmod{\mathfrak{P}_k},
\]
mapping the $1$-dimensional subspace $q_k\mathfrak{p}_k^{-1}/q_k\OO_K$ onto the $i(k)$-th coordinate line $\F_{q_k} e_{i(k)}$.
Set $B_K \deff \sum_{k=1}^\infty \mathfrak{p}_k^{-1} \subseteq K$.

Because the primes $q_k$ are pairwise distinct, localizing $B_K$ gives $(B_K)_{(p)} = (\OO_K)_{(p)} \cong \z_{(p)}^n$ for $p \notin \{q_k\}$ and $(B_K)_{(q_k)} = (\mathfrak{p}_k^{-1})_{(q_k)} \cong \z_{(q_k)}^n$ (as $\mathfrak{p}_k^{-1}$ is free of rank $n$ over $\z$), so $B_K$ is locally free.
Since $q_k\OO_K = \prod_{j=1}^n \mathfrak{p}_{k,j}$ with $n\ge 2$ distinct factors, $q_k^{-1}\notin\mathfrak{p}_k^{-1}$ and hence $1/q_k \notin B_K$.
For any $0\ne x\in B_K$, choose $m\in\N$ with $mx\in\OO_K$; if $x\in q_k B_K$ with $q_k\nmid m$, then $mx \in q_k(B_K)_{(q_k)}\cap\OO_K = q_k\mathfrak{p}_k^{-1} = \prod_{j\ne i(k)}\mathfrak{p}_{k,j}$, so $q_k^{n-1} \mid N_{K/\q}(mx) = m^n N_{K/\q}(x)$ and hence $q_k \mid N_{K/\q}(x)$, while if $x\in p B_K$ for $p\notin\{q_k\}$ with $p\nmid m$, then $(B_K)_{(p)}=(\OO_K)_{(p)}$ gives $mx \in p(\OO_K)_{(p)}\cap\OO_K = p\OO_K$, so $p^n \mid N_{K/\q}(mx)$ and $p \mid N_{K/\q}(x)$.
Since $N_{K/\q}(x)\in\q^\times$ has finitely many prime divisors, $x$ is divisible in $B_K$ by at most finitely many primes, so $B_K$ is $\bt(\z)$-homogeneous.

To prove $\End(B_K)=\OO_K$, let $\varphi\in\End_{\z}(B_K)$ and extend $\varphi$ to $\tilde{\varphi}\in\End_{\q}(K)$.
Since $\OO_K$ is finitely generated over $\z$, there exists $M\in\N$ such that $\psi \deff M\tilde{\varphi}$ satisfies $\psi(\OO_K)\subseteq\OO_K$.
Under $\Phi\colon L\otimes_{\q} K \xrightarrow{\sim} L^n$, the $L$-linear extension $\id_L\otimes\psi$ corresponds to a matrix $A_\psi = (a_{ji})_{1\le j,i\le n} \in M_n(L)$ satisfying $\sigma_j(\psi(b)) = \sum_{i=1}^n a_{ji}\sigma_i(b)$ for all $b\in K$.
For every $k\ge 1$ with $q_k\nmid M$ and $\mathfrak{P}_k$ coprime to the denominators of all entries $a_{ji}$, $\psi(B_K) = M\tilde{\varphi}(B_K) \subseteq M B_K \subseteq B_K$ implies that $\psi$ preserves both $(\OO_K)_{(q_k)}$ and $(B_K)_{(q_k)} = (\mathfrak{p}_k^{-1})_{(q_k)}$.
Multiplying by $q_k$, the induced $\F_{q_k}$-linear map on $\OO_K/q_k\OO_K$ preserves $q_k\mathfrak{p}_k^{-1}/q_k\OO_K$, so under $\pi_k$ the reduced matrix $\overline{A}_\psi = (\overline{a}_{ji}\bmod\mathfrak{P}_k) \in M_n(\F_{q_k})$ preserves $\F_{q_k}e_{i(k)}$.
Thus $a_{j,i(k)} \in \mathfrak{P}_k$ for all $j\ne i(k)$.
Since each $i\in\{1,\dots,n\}$ occurs as $i(k)$ for infinitely many $k$, each $a_{ji}$ ($j\ne i$) lies in infinitely many prime ideals $\mathfrak{P}_k$ of $\OO_L$, forcing $a_{ji}=0$.
Hence $A_\psi = \operatorname{diag}(a_{11},\dots,a_{nn})$ commutes with multiplication by $t$ ($\operatorname{diag}(\sigma_1(t),\dots,\sigma_n(t))$), so $\psi$ is $K$-linear on $K=\q(t)$ and $\tilde{\varphi}$ is multiplication by $\alpha \deff \tilde{\varphi}(1) \in K$.
Because $\alpha$ preserves $(\OO_K)_{(p)}$ for $p\notin\{q_k\}$ and $\End_{\OO_K}((\mathfrak{p}_k^{-1})_{(q_k)}) = (\OO_K)_{(q_k)}$ for all $k$, we obtain $\alpha \in \bigcap_{p\in\p}(\OO_K)_{(p)} = \OO_K$.
Conversely, $\OO_K\cdot\mathfrak{p}_k^{-1}=\mathfrak{p}_k^{-1}$ gives $\End(B_K)=\OO_K$.
Consequently, $\q\End(B_K)\cong\q\otimes_{\z}\OO_K\cong K$ is a field \textup{(}hence local\textup{)}, so $B_K$ is strongly indecomposable by \cite[Corollary 7.8]{Arnold}.
In particular, $B_K$ is indecomposable of rank $n\ge 2$ and therefore has no free summands.

Finally, since the primes $q_k$ are pairwise distinct, $T_K = B_K/\OO_K \cong \bigoplus_{k=1}^\infty (\mathfrak{p}_k^{-1}/\OO_K) \cong \bigoplus_{k=1}^\infty (\z/q_k\z)$ has cyclic $p$-components for every $p\in\p$.
Pontryagin duality applied to $0\ra\OO_K\ra B_K\ra T_K\ra 0$ gives the proper short exact sequence $0\ra T_K^\vee\ra B_K^\vee\ra\OO_K^\vee\ra 0$.
After choosing a $\z$-basis of $\OO_K$, one has $\OO_K^\vee\cong\frac{\R\cdot K}{\OO_K}\cong\left(\frac{\R}{\z}\right)^n$, so $T_K^\vee$ is a delta-subgroup of $B_K^\vee$.
By \cite[Theorem 1(3), $Y=\OO_K$]{Lewis2019}, there is a diagonal embedding $\OO_K\hookrightarrow T_K^\vee\times\R\cdot K$ whose projection into $T_K^\vee$ has dense image and for which $B_K^\vee\cong\frac{T_K^\vee\times\R\cdot K}{\OO_K}$.
The chosen $\z$-basis of $\OO_K$ induces identifications $K\cong\q^n$ and $\R\cdot K\cong\R^n$, and transporting the adelic exponential sequence along these identifications gives $B_K^\vee\cong\frac{\nala(B_K^\vee)\times\R\cdot K}{K}$.
\end{proof}

Theorem \ref{zassenhaus-realisation} associates with $K$ an irreducible adelic torus $B_K^\vee$ of archimedean dimension $[K\colon\q]$, whose nonarchimedean structure is encoded by the {\it Richman type} of $B_K$ \cite[p.12, Proposition 1.9]{Arnold}.

\begin{definition}
For a number field $K$, the \emph{nonarchimedean dimension of $K$} is defined as:
\[
\dimn K \deff \min \{ \dimn(B_K^\vee) \colon B_K \text{ is a group of integers of } K \}.
\]
\end{definition}

\begin{thm}[Nonarchimedean Line through a Number Field]\label{nonarchimedean-line-number-field}
If $K$ is a number field, then $\dimn K=1$.
\end{thm}
\begin{proof}
Let $K/\q$ be a number field.
For any group of integers $B$ of $K$, $B$ has no free summands and is therefore not finitely generated over $\z$; hence $B/\OO_K$ is infinite, so $(B/\OO_K)^\vee$ is an infinite delta-subgroup of $B^\vee$.
Since a finite product of finite groups is finite, the standard representation of $(B/\OO_K)^\vee$ contains at least one infinite procyclic factor, whence $\dimn(B^\vee) \ge 1$.

Now choose the group of integers $B_K = \sum_{k=1}^\infty \mathfrak{p}_k^{-1}$ constructed in Theorem \ref{zassenhaus-realisation}, and set $G_K\deff B_K^\vee$ and $T_K\deff B_K/\OO_K$.
Since $T_K \cong \bigoplus_{k=1}^\infty (\z/q_k\z)$ for pairwise distinct primes $q_k$, Pontryagin duality yields the delta-subgroup $T_K^\vee \cong \prod_{k=1}^\infty (\z/q_k\z)$ of $G_K$, which is an infinite procyclic group (having standard representation with a single infinite procyclic factor).
Thus $\dimn(B_K^\vee) = \dimn(T_K^\vee) = 1$, and consequently $\dimn K = 1$.
\end{proof}

\begin{example}
We close with a Pontryagin-duality computation illustrating that the dual of a familiar locally compact group need not lie wholly on either side of the usual discrete–compact correspondence.
We show $(\A/\z)^\vee \cong \zhat \times \q$, the product of a compact group and a discrete torsion-free group.
We have $\frac{\A}{\zhat+\q}\cong \frac{(\A/\zhat)}{\q}\cong\frac{(\q/\z)\times\R}{\q}\cong\frac{\R}{\z}$ so there is a proper short exact sequence $0\ra\zhat+\q\ra\A\ra\frac{\R}{\z}\ra 0$ in $\mathbf{LCA}$ with Pontryagin dual proper short exact sequence $0\ra \z\ra\A\ra(\zhat+\q)^\vee\ra 0$ in $\mathbf{LCA}$, so $\frac{\A}{\z}\cong (\zhat+\q)^\vee$; equivalently, $\left(\frac{\A}{\z}\right)^\vee\cong \zhat \times \q$.
\end{example}

\end{document}